\documentclass[12pt]{article}
\title{Bounds on Volume Increase Under Dehn Drilling Operations}
\author{Martin Bridgeman
\thanks{Research at MSRI is supported in part by NSF grant DMS-9022140.}}
\input{psfig}

\newtheorem{thm}{Theorem}
\newtheorem{defn}{Definition:}
\newtheorem{lemma}{Lemma:}
\newtheorem{conj}{Conjecture:}

\input{psfig}
\usepackage{amstex,amssymb}
 
\def\pf{{\bf Proof :}}
\def\eproof{$\blacksquare$}

\def\hypoly{hyperbolic polyhedron }
\def\bpoly{basic polyhedron }
\def\bpolys{basic polyhedra }
\def\poly{polyhedron }
\def\hypolys{hyperbolic polyhedra }
\def\chypoly{compact hyperbolic polyhedron }
\def\chypolys{compact hyperbolic polyhedra }
\def\ihypoly{ideal hyperbolic polyhedron }

\def\sc{spherical complex }
\def\scs{spherical complexes }

\def\nbd{neighborhood }
\def\a{\alpha}
\def\b{\beta}
\def\c{\gamma}
\def\d{\delta}
\def\e{\epsilon}
\def\t{\theta}
\def\f{\phi}
\def\o{\overline}

\def\scc{simple closed curve }
\begin{document}

\maketitle

\begin{abstract}
In this paper we investigate how the volume of hyperbolic manifolds increases
under the process of removing a curve, that is, Dehn drilling. If the curve 
we remove is a geodesic we are able to show that for a certain 
family of manifolds
the volume increase is bounded above by $\pi \cdot l$ where $l$ is the length
of the geodesic drilled. Also we construct examples to show that there
is no lower bound to the volume increase in terms of a linear function of a 
positive power of length and in particular volume increase is not
bounded linearly in length.
\end{abstract}

\section{Introduction}
Let $M$ and $M'$ be three-dimensional manifolds such that $M'$ is homeomorphic
to $M - \a$, where $\a$ is a \scc in $M$. Then we say that $M'$ is obtained
by Dehn drilling on $M$ and  $\a$ is the curve drilled.
 
If $M$ and $M'$ are complete hyperbolic manifolds of finite volume then 
Gromov (\cite{Th}) showed that the  volume of $M'$ is greater 
than the  volume of $M$. Therefore hyperbolic Dehn drilling always 
increases volume, but  what this increase depends on is unknown. In the
following we study the case of drilling out a geodesic and show how the 
volume increase is related to the length of 
the geodesic drilled.

\begin{conj}
Dehn drilling a geodesic increases volume by at most $\pi \cdot L$ where $L$
is the length of the geodesic to be drilled.
\end{conj}

In this paper we will prove this conjecture for a class of Dehn drillings as 
well as give motivation for why the conjecture should be true in general.
We also show that in general there is no linear lower bound on the volume 
increase by showing that there is no linear lower bound for the 
set of Dehn drillings we consider.

\section{Motivation}
The motivation for the conjecture and the idea behind the proof in the special
case to be considered is the following. If $M'$ is obtained by hyperbolic 
Dehn drilling a
geodesic $\a$ in  $M$ then we first construct a continuous one-parameter family
of hyperbolic cone-manifolds $M_\t, 0 < \t \leq 2\pi$ interpolating
between $M$ and $M'$. By interpolating we mean that $M_{2\pi} = M$, 
$\lim_{\t \to 0} M_{\t} = M'$ and $M_\t$ is homeomorphic to $M$
with  cone axis $\a$ and cone angle $\t$.

Having constructed a one-parameter family of cone-manifolds 
interpolating between
$M$ and $M'$, we apply results
of Craig Hodgson (\cite{Ho}) which describe the derivative of volume of 
a smooth 
one-parameter family of cone manifolds.

Before we give Hodgson's formula, we state the Sclafli formula for 
the variation of volume
of a one-parameter family of polyhedra in a space of constant curvature, 
on which Hodgson's formula is based.

\begin{thm}
(Sclafli, Hodgson \cite{Ho}) Let $X_t$ be a smooth one-parameter family of
polyhedra in a simply connected $n$-dimensional space of constant curvature $K$.
Then the derivative of the volume of $X_t$ satisfies the equation:
\begin{eqnarray*}
       (n-1) K dVol(X_t) & = & \sum_{F} V_{n-2}(F) d\f_F 
\end{eqnarray*}
where the sum is over all co-dimension two faces of $X_t$, $V_{n-2}$ is 
the $(n-2)$-dimensional volume and $\f_F$ is the dihedral angle at $F$.
\end{thm}

In the cone manifold case we have the following.

\begin{thm}
(Hodgson\cite{Ho}) Let $C_t$ be a smooth family of (curvature $K$) 
cone-manifold 
structures on a manifold with fixed topological type of singular locus. 
Then the 
derivative of volume of $C_t$ satisfies
\begin{eqnarray*}
       (n-1) K dVol(C_t) & = & \sum_{\Sigma} V_{n-2}(\Sigma) d\f_{\Sigma} 
\end{eqnarray*}
where the sum is over all components $\Sigma$ of the singular locus of $C$ and 
$\t_{\Sigma}$ is the cone angle along $\Sigma$.
\end{thm}
Therefore in our case where the singular locus is  a single cone axis 
$\a$ we have

\begin{eqnarray*}
       & dV  = & - \frac{1}{2} l(\t)d\t \nonumber \\
       & \Delta V = & \frac{1}{2}\int_{0}^{2\pi} l(\t)d\t \label{inteq}
\end{eqnarray*}
where $\t$ is the cone  angle along $\a$, $l(\t)$ is the
length of $\a$ in $M_{\t}$ and where $\Delta V$ is the 
volume increase under Dehn drilling.

Therefore if the function $l$ is monotonically increasing in $\t$
then it follows that
$$ \Delta V \leq \pi\cdot L$$
where $L$ is the length of the geodesic $\a$ in $M$.

The conjecture and the proof of the conjecture for the class of Dehn 
drillings we consider, is based on this two-fold approach of interpolation
and monotonicity.

\section*{Basic Polyhedra}
\begin{figure}
\begin{center}~
\psfig{file=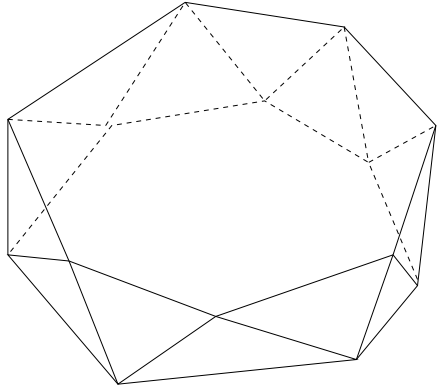,height=2in}
\end{center}
\caption{Ideal hyperbolic drum \label{drum}}
\end{figure}

We will now prove that for a Dehn drilling operation on {\em basic polyhedra},
the conjecture is true. This Dehn drilling operation has some interesting
properties including giving a new way to enumerate all \bpolys and hence all
link projections (\cite{B2}).
A  basic polyhedron is an \ihypoly $P$ with all dihedral angles  
being right angles. Therefore all vertices have valence four. 
Let $\mathcal{B}$ be the set of all basic polyhedra. 
Examples of basic polyhedra are the ideal hyperbolic drums
$\{T_n\}_{n \geq 3}$ which have two faces being regular
ideal $n$-gons, corresponding to the top and bottom of the drum 
and every other face being an ideal triangle (figure \ref{drum}).

As all the dihedral angles are  $\pi/2$,  basic polyhedra can 
also be thought of as  orbifolds.
In \cite{B1} it was shown that if $P$ is a basic polyhedron then it has a
four-fold cover (in the orbifold sense) $L_P$ which is a hyperbolic link 
complement in $S^3$ having a boundary component for each vertex of $P$.

\subsection*{Surgery and Dehn drilling}
\begin{figure}
\begin{center}~
\psfig{file=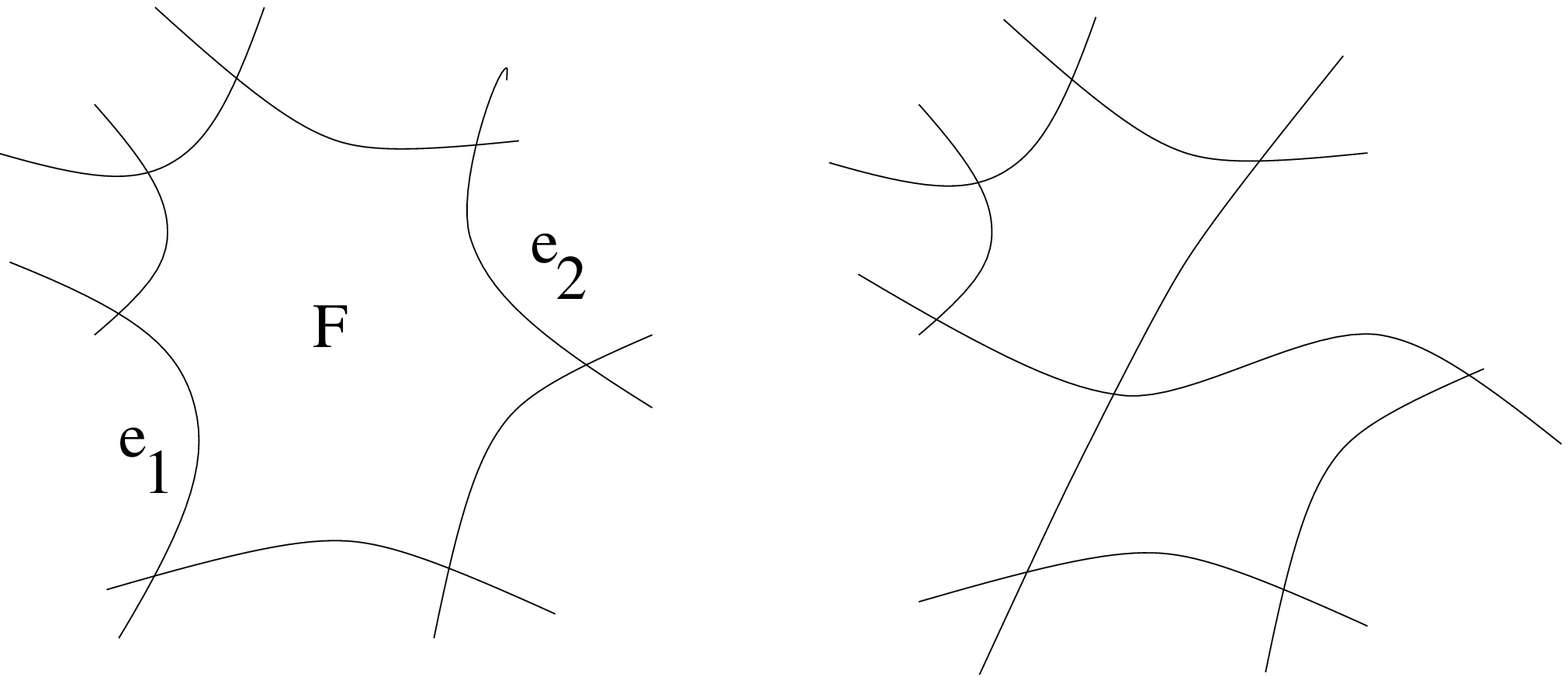,height=2in}
\end{center}
\caption{Surgery \label{surgery}}
\end{figure}

If $P$ is a basic polyhedron and $F$ is a non-triangular face of $P$ then
we define {\em surgery} on $P$ as follows. Choose two non-adjacent
edges $e_1,e_2$ of $F$ and pinch them together to form a new combinatorial
\poly $\o{P}$ (figure \ref{surgery}). By Andreev's theorem (\cite{Th}), $\o{P}$
has a realization $P'$ which is an \ihypoly with all dihedral angles right
angles and therefore is also a basic polyhedron.
We say $P'$ is obtained by surgery on $P$ and therefore surgery is an operation
defined on the set $\mathcal{B}$.

In the sense of orbifolds, $P'$ is obtained by Dehn drilling the unique 
geodesic perpendicular to both $e_1$ and $e_2$ in $P$ and if we lift to the
four-fold covers, we see that
$L_{P'}$ is obtained by Dehn drilling  a simple closed geodesic 
in $L_{P}$ (\cite{B2}). Thus surgery corresponds to a Dehn drilling operation.

We define a partial ordering on $\mathcal{B}$ by saying 
$P \prec P'$ if $P'$ is the
result of repeated surgeries starting with $P$. This partial ordering has 
as initial objects the set $\{T_n\}_{n \geq 3}$ mentioned above (\cite{B2}) 
and therefore we can enumerate all basic polyhedra by repeated surgery 
on these initial objects.

\begin{thm} 
Surgery increases volume by at most $\pi \cdot L$ where $L$ is the
length of the geodesic drilled.
\label{mainthm}
\end{thm}

To prove this theorem we must first interpolate between $P$ and $P'$ where
$P'$ is obtained by surgery on $P$.

\section{Interpolating Surgery}
Let $P$ be a \bpoly with $F$ a non-triangular face of $P$ and $e_1,e_2$ be two
non-adjacent edges of $F$. Also let $g$  be the unique geodesic in $F$ 
perpendicular to $e_1,e_2$. 
To interpolate between $P$ and $P'$  we need a continuous one-parameter family 
of polyhedra $P_\t, 0 < \t \leq \pi$ 
which have combinatorial type of $P \cup \{g\}$, with dihedral angle $\t$ 
along $g$ and all other edges having the same dihedral angles
as the corresponding edge of $P$ and
with $P_{\pi}= P, \lim_{\t \to 0}{P_{\t}} = P'$.
Intuitively $P_\t$ is obtained by bending $P$ along $g$ an amount $\pi - \t$. 
To find such a one-parameter family of polyhedra we use results
of Hodgson and Rivin (\cite{HR}) which characterizes compact hyperbolic 
polyhedra in terms of their {\em spherical duals}.

\subsection*{Spherical Dual of a Polyhedron}
\begin{figure}
\begin{center}~
\psfig{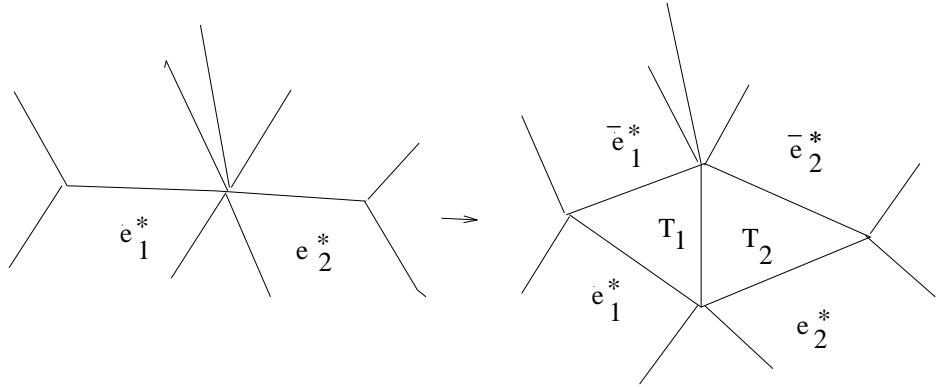}
\end{center}
\caption{Cutting open $S$ to make $S_\t$ \label{bigon}}
\end{figure}

If $P$ is a \hypoly then we form its spherical dual $P^*$ which is 
a spherical complex combinatorially dual to the polyhedron $P$ as follows. 
Each vertex
corresponds to a spherical polygon with edges having lengths being the
exterior dihedral angles of the edges incident to the vertex and with face
angles being the complementary angles of the face angles at the 
vertex (complement of $\a = \pi - \a$).
Similarly an ideal vertex corresponds to a spherical hemisphere with vertices 
dividing the equator into arcs of length less than $\pi$.
Let $S = P^*$ and define  \scs $S_\t$ by
cutting $S$ open at $e^*_1,e^*_2$ (dual edges to $e_1,e_2$), and attaching
a spherical bigon $B$ with angle $\pi - \t$. $B$ is further divided into two
triangles $T_1,T_2$ by joining the two points on the equator by an edge
$g^*$. If $P_{\t}$ is realisable as a hyperbolic polyhedron then its
spherical dual is $S_{\t}$  with
$g^*$ dual to $g$. Note that each of $e^*_1,e^*_2$ split into two edges 
$e^*_i,\o{e}^*_i ,i=1,2$ (figure \ref{bigon}). Thus we have a continuous 
one-parameter 
family of \scs $S_\t, 0 \leq \t \leq \pi$ which interpolates between the 
spherical duals of $P$ and $P'$. Before we can interpolate between $P$ and 
$P'$ we need to know more about the map between \hypolys and spherical 
complexes.

The following theorem due to Craig Hodgson and Igor Rivin gives a
characterization of \chypolys in terms of their spherical duals.
\begin{thm}
(\cite{HR}) A metric space ($M,g$) homeomorphic to $S^2$ can arise as the
spherical dual $P^*$ of a \chypoly $P$ in $H^3$ if an only if the following 
hold.

(a) The metric has constant curvature 1 away from a finite collection of
	cone points $c_i$.

(b) The cone angles at $c_i$ are greater than $2\pi$.

(c) The lengths of closed geodesics of ($M,g$) are all strictly greater than
    $2\pi$.
\label{char}

Moreover $P$ is uniquely determined.
\end{thm}

\subsection*{Approximating by compact Polyhedra}
Since the polyhedron $P_{\t}$ is not compact, 
we will find a family of \chypolys  $P^t_\t, t>0$ approximating $P_{\t}$ by
applying the above result.
We define \scs $S^t_\t$ by modifying $S_\t$ so that it is the dual
of a \chypoly as follows. We increase the length of every edge of $S_\t$ which
does not belong to $B$ by a factor of ($1+t$) and add a vertex $v_i$ at the 
pole
of each hemisphere $H_i$, and join $v_i$ by an edge to each vertex on 
the equator of $H_i$. Note that if we let $t=0$ then 
$S^{0}_{\t} = S_{\t}$. It will be shown, using the above 
characterization (theorem \ref{char}), that $S^{t}_{\t}$ is the dual
of a \chypoly $P^t_\t$. Furthermore $P^t_\t$ has combinatorial type of 
$P_\t$ with its ideal
vertices truncated and the truncating faces 
are perpendicular to the faces of $P^t_\t$ intersecting it. Thus in the Klein 
model, all the planes 
intersecting a truncating face meet in a single point outside
$S^2_{\infty}$ corresponding to the unique 
point perpendicular to the truncating face. A 
{\em generalized hyperbolic polyhedron} is when we allow vertices to be 
outside $S^2_{\infty}$ in the Klein model and we call the vertices outside
of $S^2_{\infty}$ {\em hyper-ideal}. In this paper we will not 
distinguish 
between a generalized polyhedron and the unique polyhedron obtained by
truncating the hyper-ideal vertices.
Thus $P^t_\t$ corresponds to a generalized hyperbolic polyhedron with 
combinatorial type
$P_\t$ and with the vertices corresponding to ideal vertices of
$P_\t$ being hyper-ideal. Intuitively as $t$ tends to 0 these
hyper-ideal vertices tend to $S^2_{\infty}$ and become ideal.

\begin{lemma}
Given $\t, 0 < \t \leq \pi$ then  there exists  $\d > 0$ such that
the \sc $S^t_\t$ 
is the dual of a unique \chypoly $P^t_\t$ for $0 < t < \d$.
\label{Pexists}
\end{lemma}

\pf
Obviously $S^t_\t$ is a \sc satisfying (a) and (b) of the characterization
(theorem \ref{char}),
 and therefore we need only show that there are no geodesics of 
length $\leq 2\pi$.
Let $\c^t$ be a geodesic in  $S^t_\t$. If $\c^t$ intersects the interior
of a hemisphere $H$ in a single arc then $l(\c^t |_{H^o}) = \pi$.
Therefore if $\c^t$ intersects more than two hemispheres in their interiors
then $l(\c^t) \geq 3\pi$.

{\bf Case 1: $\c^t$ doesn't intersect the interior of bigon $B$.} Then we 
get a curve in $S^t_{\pi}$ by removing the bigon and stitching back the edges.
This curve we also call $\c^t$ and it is a geodesic in $S^t_{\pi}$.

If $\c^t$  intersects the interior of two hemispheres $H_1,H_2$ and if 
$\c^t$ doesn't contain an edge then $H_1,H_2$  intersect 
in points a distance $\geq \pi$ apart giving a contradiction. Thus
$\c^t$ must contain an edge and $l(\c^t) \geq 2\pi + \pi/2 > 2\pi$.

If $\c^t$ intersects the interior of only one hemisphere $H$, then 
it must enter  and leave $H$ along an edge. We replace the arc in $H^o$
by  an edge path of the same length by going along the boundary of $H$. This
gives us a new geodesic which we still call $\c^t$ satisfying
\begin{eqnarray*}
l(\c^t) & \geq  & \pi/2 + \pi/2 + \pi/2(1+t) + \pi/2(1+t) = 2\pi + t\pi > 2\pi
\end{eqnarray*}
Therefore all geodesics in $S^t_\t$ not intersecting the interior of $B$
have length bounded away from $2\pi$.

{\bf Case 2:  $\c^t$ intersects interior of $B$.} There are a number of 
possibilities.

If $\c^t$ enters $B$ at either the apex of $T_1$ or $T_2$ then it has
to leave at the other apex and $l(\c^t |_{B^o}) = \pi$. We replace 
this arc by edge path $e^*_1 \cup e^*_2$ or $\o{e}^*_1 \cup \o{e}^*_2$
depending on which gives a geodesic as a 
result (at least one must give a geodesic). Thus we have a geodesic
which has the same length as $\c^t$ and doesn't intersect the interior of 
$B$. Therefore $\c^t$ has length $ > 2\pi$ by the first case.

If $\c^t$ intersects the interior of two hemispheres $H_1,H_2$ and both 
are  not incident to $B$ then $\c^t$ must contain an edge of 
$S^t_\t$ and therefore $l(\c^t)  \geq 2\pi + \pi/2  > 2\pi$. 
If either $H_1,H_2$ are
incident to $B$ only at an apex then $\c^t$ intersects $B$ in an arc from
one apex to the other and has been already shown to be $> 2\pi$.
Therefore  either $\c^t$ contains an edge and is longer than
$2\pi$ or $H_1,H_2$ intersect at a point $p$ along their boundaries where
$\c^t$ traverses from $H_1$ to  $H_2$. But $H_1,H_2$ correspond to 
hemispheres in $S^t_{\pi}$  whose intersection is contained in 
$e^*_1 \cup e^*_2$
and thus $H_1,H_2$ cannot intersect at point $p$. Therefore any 
geodesic in $S^t_\t$ intersecting
 the interior of two or more hemispheres has length bounded away from $2\pi$.

\begin{figure}
\begin{center}~
\psfig{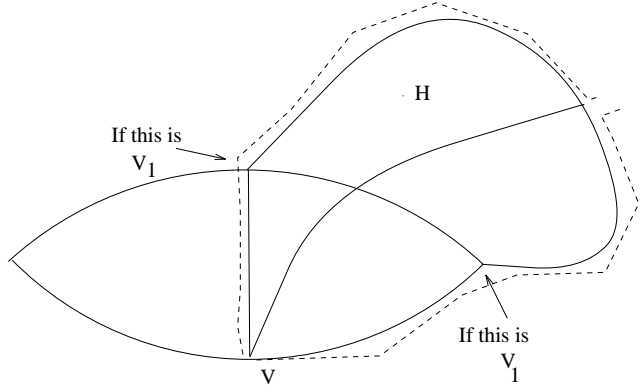}
\end{center}
\caption{Replacing geodesic by edge geodesic \label{sc}}
\end{figure}

If $\c^t$ only intersects the interior of one hemisphere $H$ and  it
enters $B^o$ along an edge and leaves along an edge then it must intersect 
$B^o$
along the equatorial edge. Therefore there is a geodesic in $S^t_{\pi}$ of 
length
$l(\c^t) - (\pi - \t)$ obtained by removing the bigon as before and 
taking what is left
of $\c^t$ as the curve. This must be of length $ > 2\pi$, therefore $\c^t$ has
length $\geq 2\pi + (\pi -  \t) > 2\pi$. If $\c^t$ enters through $H$ 
and leaves along
an edge at vertex $v$ of $B$ ($v$ is not an apex), then we first
replace the interior arc in $H$ by a boundary edge path $b$ of $H$ of the
same length, we call this new curve $\c^t$ also.
Now $\c^t$ contains two sides of a spherical triangle $\o{T}$ in $B$. $H$ 
intersects $B$ in an edge $e$ and let $v_1$ be the vertex of $B$ contained in 
$b \cap e$. Now replace the part of $\c^t$ consisting of 
two sides of the spherical triangle $\o{T}$
by the edge of $B$ joining $v,v_1$ (figure \ref{sc}). This gives a geodesic of
smaller length that is an edge path and therefore $l(\c^t) > 2\pi$.

If $\c^t$ intersects no hemisphere in its interior then $\c^t$  is an edge
geodesic containing the equatorial edge of $B$. Removing the bigon
we get a geodesic in
$S^t_{\pi}$ of smaller length which implies $l(\c^t) > 2\pi$.

Since no geodesic of $S^t_\t$ has length $\leq 2\pi$,
then 
$S^t_\t$ is  the spherical dual of a unique 
\chypoly $P^t_\t$ (by theorem \ref{char}). \eproof \newline

To find what is the combinatorial structure of $P^t_\t$ we need to find a
cell division of $S^t_\t$ into convex spherical polygons by geodesics of length
less than $\pi$ and whose endpoints are cone points. This is how we described
$S^t_\t$ originally, but this may not be the only cell division possible. To
show that it is, consider the apex $v$ of one of the triangles $T_1,T_2$. The 
only other vertices  a distance less than $\pi$ from $v$ are those
vertices connected to $v$ by an edge in the original cell division and 
therefore the edges from $v$ 
must be a subset of the original edges from $v$. Also all edges must be 
included as if any are left out
then there is an angle $\geq \pi$ subtended at $v$. Now considering a vertex 
$v_1$ on the equator we have an angle of $\pi$ inside the bigon $B$ at $v_1$. 
Therefore the equatorial edge of $B$ must also be included in the division.
Using the fact that any geodesic intersecting the interior of a  hemisphere
must have length $\geq \pi$ and that no vertex can have an angle of $\geq \pi$
subtended at it, it is easy to show that all other edges of the original
division are necessary. There cannot be any more edges as any other edge from
vertex to vertex would cross another edge in an interior point which can't
happen. Therefore $P^t_\t$  has dual given by the original cell division of
$S^t_\t$.

$P^t_\t$ is  a general \hypoly with  
combinatorial type $P \cup g$ and  hyper-ideal vertices corresponding
to the ideal vertices of $P$. The dihedral angles between a
truncating plane and the planes it intersects are all $\pi/2$ . Also the 
dihedral angle along the edge corresponding to $g$ is $\t$ and except 
for the
four edges incident to $g$ which have dihedral angles $\pi/2$, all others have
dihedral angle $\pi/2(1-t)$. Similarly the face angles on all truncating faces
are $\pi/2(1-t)$ except where the four edges incident to $g$ intersect. 
The face angles on the two faces perpendicular to $g$, at the vertices 
incident to $g$,  are $\t$, and all other face angles are $\pi/2$.
As $t$ tends to zero all  dihedral angles
approach $\pi/2$ except along $g$ which has fixed dihedral angle $\t$.

\section{Monotonicity and Volume Increase}
Lemma \ref{Pexists} shows that $P^t_{\t}$ exists for $0 < t \leq k(\t)$,
$0 < \t \leq \pi$ where $k(\t) > 0$. We can choose $k(\t) = c \cdot \t$
where $c$ is a positive constant.
Hodgson and Rivin showed (\cite{HR}) that the map from the space of
\chypolys to their spherical duals is a homeomorphism and hence $P^t_{\t}$
is a continuous two-parameter family of polyhedra in the above domain. 
When all the dihedral
angles are acute we have continuity by Alexandrof (\cite{A}) and therefore
$lim_{t \to 0}{P^t_{\pi}} = P$.  Thus we extend $P^t_{\t}$
by letting $P^0_{\pi} = P$. Also if $0 < \t \leq \pi/2$   
then $lim_{t \to 0}{P^t_{\t}} = P^0_{\t}$ where $P^0_{\t}$ is the 
unique \hypoly with combinatorial type of $P \cup \{g\}$ and dihedral
angles $\pi/2$ except along $g$ which has  dihedral angle $\t$. We can
now take the limit as $(t,\t)$ tends to $ (0,0)$ and we get that
$lim_{(t,\t) \to (0,0)}{P^t_{\t}} = P'$.
Therefore $P^t_{\t}$ is well defined and continuous on the triangle 
$0 \leq t \leq  c \cdot \t$, $0 \leq \t \leq \pi$ except for along the open
interval of length 
$\pi/2$ from $(0,\pi/2)$ to $(0,\pi)$.

\begin{figure}
\begin{center}~
\psfig{file=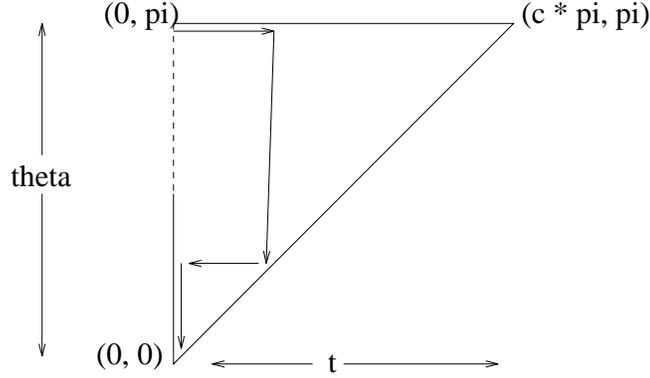,height=2in}
\end{center}
\caption{2-parameter domain and path from $P$ to $P'$ \label{1par}}
\end{figure}

We now choose a continuous one-parameter family of polyhedra interpolating
between $P$ and $P'$ by taking the image of a piecewise linear path joining
$(0,\pi)$ and $(0,0)$ in the domain of $P^{(\cdot)}_{(\cdot)}$ 
(figure \ref{1par}). Thus the one-parameter family $X$ interpolating between
$P$ and $P'$ is the product of four continuous one-parameter families 
$X_{i}, i=1,\ldots,4$ obtained by taking the image of the path
joining the four points 
$(0,\pi),(0,t),(c/t,t),(c/t,0),(0,0)$ by alternately horizontal and vertical 
lines. 

The volume change from $P$ to $P'$ $\Delta V$  satisfies 
$$\Delta V = \Delta V_1 + \Delta V_2 + \Delta V_3 + \Delta V_4$$
where $\Delta V_i$ is the volume change over $X_i$.

By continuity given $\epsilon > 0$ we can choose a $t > 0$ such that
\begin{eqnarray*}
\vert \Delta V_i \vert & \leq  \epsilon & i \neq 2 \\
\vert l^t(\pi) - l^0(\pi) \vert & \leq & \epsilon
\end{eqnarray*}
where $l^0(\pi) = l$, the length of the original geodesic in $P$.

If $l^t$ is monotonic then
\begin{eqnarray*}
0 \leq \Delta V_2 & \leq & 1/2(\pi - t/c)l^t(\pi) \leq \pi/2 \cdot l^t(\pi) \\
\Rightarrow 0 \leq \Delta V & \leq & \pi/2 \cdot l^t(\pi) + 3\epsilon \leq 
\pi/2 \cdot l + \tilde{\epsilon}
\end{eqnarray*}

Therefore by choosing smaller values of $t$ we have 
$$\Delta V \leq \pi/2 \cdot l$$
Taking the four-fold cover of $P$ given by the hyperbolic link complement $L_P$
we have  
$$\Delta V' = 4 \cdot \Delta V \leq 4\cdot(\pi/2 \cdot l) = 
2\pi l=\pi L$$
where $\Delta V'$ is the volume increase
from  $L_P$ to $L_{P'}$
and $L$ is the length of the geodesic drilled from $L_P$ to obtain $L_{P'}$.

\subsection*{Monotonicity}
The following  theorem is based on Cauchy Rigidity and is the step
needed to prove monotonicity, but first a  
definition.
\begin{defn} A {\em signed polyhedron} $P$ is a polyhedron whose edges are
labeled with either  $+$ , $-$ or $0$. 
\end{defn}

A vertex $v$ of a signed polyhedron has
{\em type} given by listing the signs on the edges incident to 
$v$ in a clockwise
or anti-clockwise manner (with the usual equivalence) and a face $F$ is called
a zero-face if all its edges are marked with zero.

\begin{thm}
Let $P$ be  a trivalent signed polyhedron such that, except for two 
disjoint faces 
$F_1,F_2$, every face is either a zero-face or has at least four sign changes.
Then the space obtained by collapsing all zero-edges $\tilde P$ 
is topologically a bouquet of spheres
where either

(1) Each sphere has the trivial cell division (a point and a disc). 

or 

(2) There is exactly one sphere $\tilde S$ in the bouquet with
a non-trivial cell division. $\tilde S$ has two cells with no sign
changes and every other cell has exactly four sign changes. Also all 
vertices of $\tilde S$ have 
type $(++-)$ , $(--+)$  or  $(+-+-)$ and therefore valence three or four.
\label{cauchy}
\end{thm}

\begin{figure}
\begin{center}~
\psfig{file=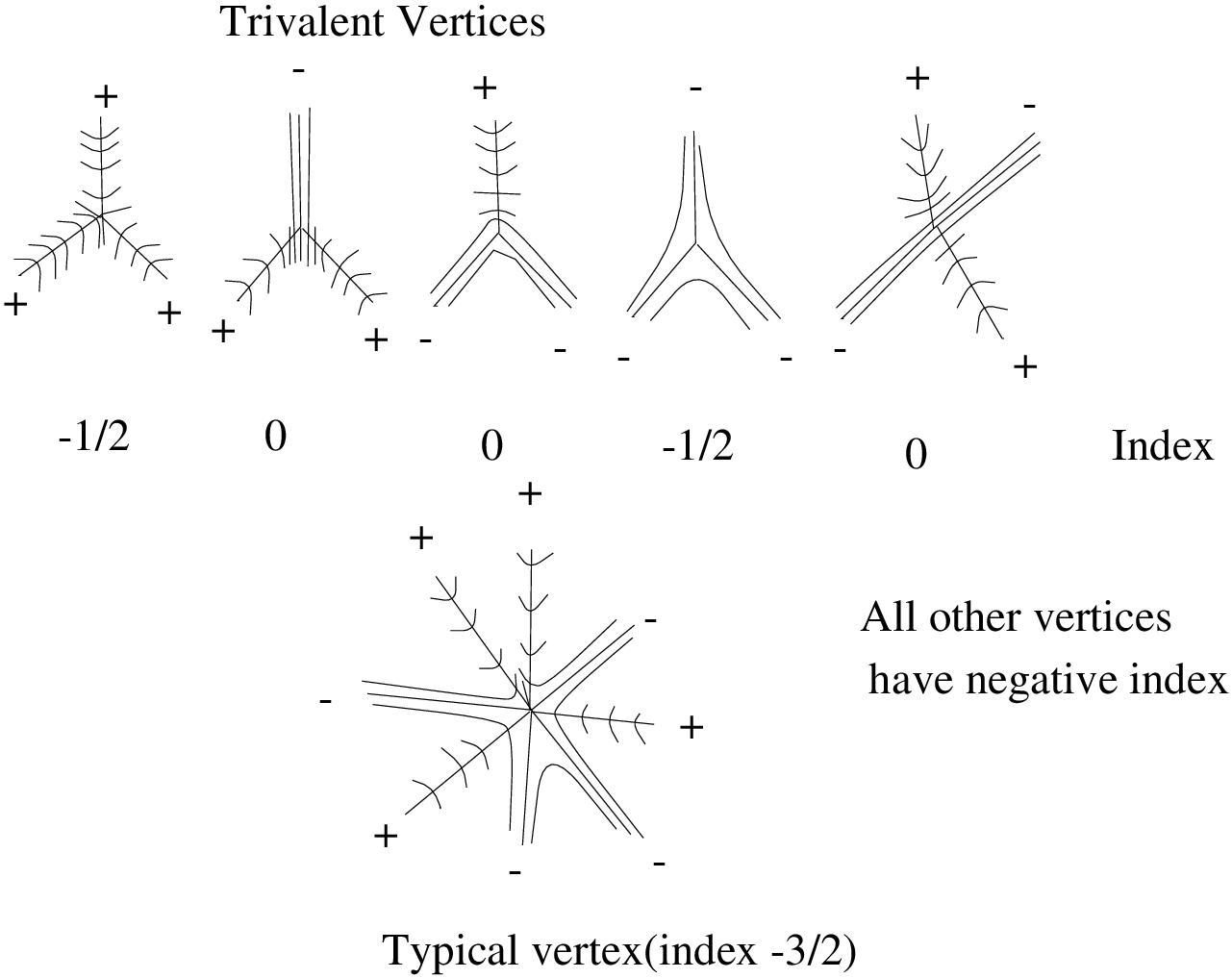,height=3in}
\end{center}
\caption{Indices of vertices\label{vindex}} 
\end{figure}

\pf Let $\Sigma$ be the set of edges of $P$ and 
$\Sigma = \Sigma_+ \cup \Sigma_- \cup \Sigma_0$ be the partition into $+$ , $-$
and  $0$ edges.
If $\vert \Sigma_0 \vert = 0$ then construct  a vector field on  $P$
as follows. In a \nbd 
of $\Sigma$ the vector field is transverse to $\Sigma_+$ and parallel 
to $\Sigma_-$. The index of a vertex $v$ can easily be calculated from 
its type.
One way is to read the type of a vertex from left to right, count $-1/4$ 
for sign changes, $-1/2$ for same sign  and then adding $1$. 
This is never positive if the valence is at least three and if the 
valence of $v$
is greater than four its index is always negative. In our case 
the vertices all have valence three but in general only 
vertices of type $(++-) , (--+)$ and $(+-+-)$ (which have index zero), have 
non-negative indices (figure \ref{vindex}).

\begin{figure}
\begin{center}~
\psfig{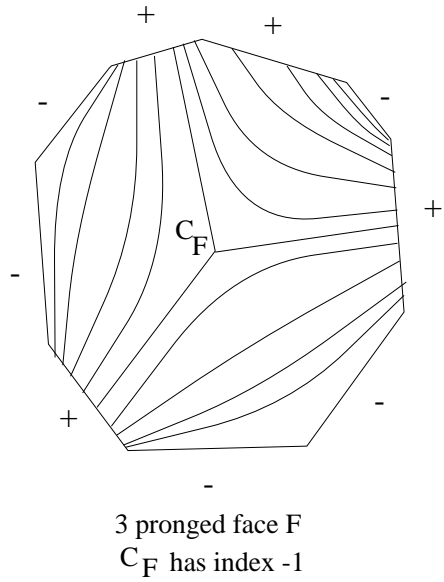}
\end{center}
\caption{Extend vector field over face \label{findex}}
\end{figure}

The vector field is extended over a  face $F$ as follows. As we go round the
edges of $F$ the sign changes $2n$ times , which means that
there are $n$ disjoint $+$ segments joined together by $n$ disjoint $-$ 
segments. Choose a center point $c_F$ in the interior of $F$ and for each 
$+$ segment we join an interior point of it to the center
$c_F$ by a prong and extend the vector field in a \nbd of these $n$ prongs in 
the obvious way. Cutting along the  prongs divides $F$ into $n$ pieces over
which the vector field can easily be extended. The index of the critical point
$c_F$ is $1 - n/2$ (figure \ref{findex}). Note that if the face has no sign 
changes then the vector field
can  still be extended with $c_F$ having index 1, in keeping with the formula.
The interior of each face has at most one critical point and all other
critical points are vertices of $P$. Each vertex  contributes a non-positive 
index and 
each face $F \neq F_1,F_2$ has at least four sign changes, which implies 
$c_F$ has  non-positive index. In order to have the sum of the indices 
equal to two (the Euler number of $P$), the indexes of both 
$c_{F_1},c_{F_2}$ must be one and therefore they are of constant sign. 
Also all vertices and centers $c_F$ , $F \neq F_1,F_2$
must have index zero. Therefore all vertices have type $(++-) , (--+)$ 
and each $F$ has exactly four sign changes.

If $\vert \Sigma_0 \vert \neq 0$  then we cannot find a vector field extending 
over $\Sigma_0$. We  avoid  this problem by collapsing these zero-edges. First
consider $\o{\Sigma}_0 \subseteq \Sigma_0$ to be the collection of edges 
belonging to loops in $\Sigma_0$. This decomposes $P$ into regions $\{R_i\}$
which have disjoint interiors and are $n_i$-connected regions whose boundaries
consist of zero-edges. Some regions correspond to a zero-face which we call
zero-regions. To find a cell division with no zero-edges we must collapse
$\Sigma_0$. Each region $R_i$ becomes a sphere $\tilde{S}_i$ in a bouquet 
of spheres $\tilde P$. The zero-regions become spheres with the trivial cell
division of a point and a disc and every other face $F$ collapses to a cell 
$\tilde F$ of some $\tilde{S}_i$.

\begin{figure}
\begin{center}~
\psfig{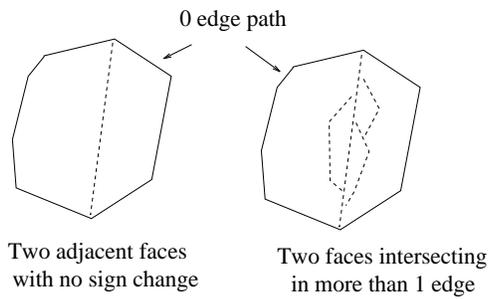}
\end{center}
\caption{Valence after collapse is at least three \label{collapse}}
\end{figure}

We collapse ${\Sigma}_0$ by first collapsing $\o{\Sigma}_0$. 
Each $R_i$ becomes a sphere $S_i$. If $R_i$ is not a zero-face then we have
each component of $\partial R_i$, $\b_{ij}$ collapses to a vertex $v_{ij}$
on a sphere $S_i$. As $R_i$ is not a face there is an edge path $\c_{ij}$
splitting $R_i$ into two regions $R^1_{ij},R^2_{ij}$ with endpoints on
$\b_{ij}$. If there are no other paths from $\b_{ij}$ to the interior then
either $R_i$ is the union of two faces or there are two faces which intersect
in two separate edges (figure \ref{collapse}). Both give contradictions. 
Therefore collapsing $\b_{ij}$ gives vertex $v_{ij}$ with valence at least
three.
Thus the vertices on $S_i$ after collapsing $\o{\Sigma}_0$, have valence 
at least three. If we now collapse $\Sigma_0 - \o{\Sigma}_0$ to get sphere 
$\tilde{S}_i$ in $\tilde{P}$ the valences cannot go down as  the only
zero-edges left on $S_i$ form a collection of disjoint trees.

Therefore if $\tilde{S}_i$ does not correspond to a zero-face it has a 
decomposition into faces all but two of which have at least four
sign changes and all 
vertices have valence at least three. Thus as before we can give $\tilde{S}_i$
a vector field which implies that there can be at most one sphere 
$\tilde S$ which
has a non-trivial cell division. On $\tilde S$ the vertices all have type 
$(++-) , (--+)$ or $(+-+-)$ , the faces corresponding to $F_1,F_2$ give 
rise to 
faces $\tilde{F}_1 , \tilde{F}_2$ of $\tilde S$ which have constant sign on 
their boundaries and all other faces have exactly four sign 
changes. \eproof \newline 

Now we apply the above theorem. If we compare the corresponding faces 
of $P^t_{\t_1}$ and
$P^t_{\t_2}$ then apart from the two faces $F_1,F_2$ which are incident and
perpendicular to $g$ , all others have identical face angles. 
Let $Q$ be the marked polyhedron combinatorially equivalent to both 
$P^t_{\t_1}$, $P^t_{\t_2}$ and marked  with $+$ , $-$ or $0$ 
depending on whether
the length of the edge in $P^t_{\t_1}$ is bigger, smaller or equal in length
to the corresponding edge of $P^t_{\t_2}$. By the below result
in planar hyperbolic
geometry  (lemma \ref{tech}), any face $F \neq F_{1}, F_{2}$ is either 
a zero-face or has at least four sign changes.
Therefore by the above (theorem \ref{cauchy}), both faces $F_1,F_2$ have no 
sign changes. If $\t_2 > \t_1$ then again using the below lemma,
both faces must have constant sign $+$. If the edge $g$ is  $+$ then
we get two $(+++)$ vertices
and if the edge $g$ is $0$ then after collapsing we get a $(++++)$ vertex, 
both giving contradictions. 

Therefore $g$ 
must be marked with $-$ and thus $l^t$ is monotonic in $\t$.

Theorem \ref{cauchy} gives more information about the pattern on a marked 
polyhedron $P$. The region $R$ corresponding to $\tilde S$ is obtained by 
removing the zero-faces of $P$ and has no zero-loops with interior edges. 
When we collapse the boundary curves of $R$ to get sphere $S$ then any vertex
of $S$ has valence three or four. Any zero-edge path of more than one edge will
give a vertex in $\tilde S$ with valence greater than four and therefore
every zero-edge in $S$ is
isolated. These isolated zero-edges in $S$ must join two trivalent vertices
which have alternating signs as you go around a \nbd of the zero-edge.

\begin{figure}
\begin{center}~
\psfig{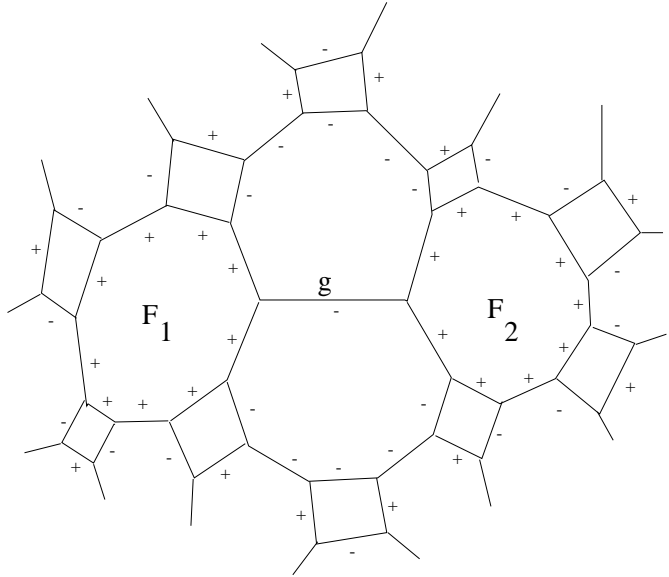}
\end{center}
\caption{How edge lengths in $P^t_{\t}$ vary with $\t$ \label{pattern}}
\end{figure}

In the case of the markings on $Q$ we obtain figure \ref{pattern}.

 The following is the technical lemma needed above. A proof can
be found in the referenced paper.

\begin{lemma} (\cite{HR}) Let $E_1,E_2,\ldots,E_n$ be lines in $H^2$ defining
a convex polygonal curve $C$. Let $A_i = E_i \cap E_{i+1} (i = 1,\ldots,n-1)$ be
the vertices of $C$ , let $\nu _i = E_i^{\perp} (i = 1,\ldots,n)$ be the outward
unit normals to  the Minkowski inner product $\langle \cdot,\cdot \rangle$, 
and let $l_i$ be the length of $E_i \cap C = A_{i-1}A_i (i = 1,\ldots,n-2)$. 
Let $E'_1,E'_2,\ldots,E'_n$ be lines in $H^2$ defining another such curve
$C'$ with $A'_i, \nu'_i,l'_i$ defined as above. Assume that

(i) $<A_i = <A'_i$ for $i = 2,\ldots,n-1$ and

(ii) $l'_i \geq l_i$ for $i = 2,\ldots,n-1$.

Then
\begin{eqnarray*}
\langle \nu'_1 , \nu'_n \rangle & \leq & \langle \nu_1 , \nu_n \rangle ,
\end{eqnarray*}
with equality if and only if $A_{1}A_{2}{\ldots}A_{n-1}$ is congruent to 
$A'_1A'_2{\ldots}A'_{n-1}$.
\label{tech}
\end{lemma}

An easy consequence of this lemma is the following. Let
$F_a,F_b$ be  two hyperbolic $n$-gons with face angles equal. If 
we label an $n$-gon $F$  with $+$, $0$ or $-$ by comparing lengths of 
corresponding sides of $F_a,F_b$, then  either $F$  has all sides 
marked with $0$, in which case $F_a,F_b$ are congruent, or, $F$ has at least 
four sign changes as we
go around the perimeter. This shows that $Q$ satisfies the 
hypothesis of the theorem. 

To show that $F_1,F_2$ are marked $+$ we need to relate the Minkowski inner 
product to face angles.
If two edges of a compact hyperbolic polygon intersect in an angle $\t < \pi$
then the inner product of their outward normals is $\cos(\t_{ext})$ where 
$\t_{ext}$ is the exterior angle of $\t$.
Now the convex polygons $p_1,p_2$ corresponding to $F_1$ in 
$P^t_{\t_1},P^t_{\t_2}$ respectively, 
have the same face angles except at the vertex incident to $g$, where
$p_i$ has the  angle  $\t_i$. Thus if the sides of 
$p_1$ where all smaller or
equal in length to the sides of $p_2$  then we would have by the above lemma
that the inner product of the outward normals of the two sides intersecting at
$g$ in $p_1$ is greater than or equal to the same inner product in $p_2$, that 
is $\cos(\pi - \t_{1}) \geq \cos(\pi - \t_{2})$. 
But $\t_1 < \t_2$,
and therefore the sides of $p_1$ must be greater than or equal in 
length to the corresponding
sides in $p_2$. This gives the needed step to show that both $F_1,F_2$ are
of constant sign $+$.

\section{Combinatorial Bound on Volume Increase}
Here we will show that bending along a geodesic $g$ on a face $F$ which is
perpendicular to two edges of $F$ increases volume at most  linearly in
the number of edges of $F$. 

The idea is that as $\t$, the exterior dihedral angle at $g$ varies, 
$F$ is still intrinsically
an ideal n-gon. If $l(\t)$ is long then this forces $F$ to be narrow
somewhere across $g$. This implies that the two faces $F_1,F_2$ obtained by
bending $F$ along $g$ (which meet with exterior dihedral angle $\t$), 
have faces $F_3,F_4$ adjacent to $F_1,F_2$ respectively which are connected 
by a short path on the bent face $F$. As $F_3,F_4$ cannot intersect,
this shows that $l(\t)$ must decrease as we bend otherwise $F_3,F_4$ would
be forced to intersect.

We will describe the proof for when we have 
a continuous one-parameter family of polyhedra $P_{\t}$ interpolating between
$P$ and $P'$. We obtain the result in general by approximating as before by
the one-parameter family $X^t$. 

\subsection*{Long Geodesics give Thin Polygons}

An ideal hyperbolic quadrilateral has two perpendiculars and it can be shown
using elementary hyperbolic geometry that the hyperbolic sine's of the half
lengths of these are reciprocal. The following lemma is a generalization of
this fact.

\begin{lemma}
If $F$ is an ideal hyperbolic n-gon and $g$ the unique perpendicular to edges
$e_1,e_2$ of $F$ then  $\exists$   $e_3,e_4$ separated by $g$ with distance
$d$ between them satisfying
$$\sinh(d/2) \cdot \sinh(\frac{(\pi - 2)l}{2\pi(n-3)}) \leq 1$$
where $l$ is the length of $g$.
\label{thin}
\end{lemma}

\pf
Choose $\e > 0$ such that the one-sided $\e$-\nbd of $g$ has the same area as $F$.$$\Rightarrow l \sinh(\e) = (n-2)\pi$$
If we take $N_{\e}$ the two-sided $\e-$\nbd of $g$ then for each 
$x \in g$ there is
an arc $\c_x$ perpendicular to $g$ centered at $x$ of radius $\e$. We define  
three subsets of $g$ as follows:
 
$$g_i  =  \{ x \in g : \vert \c_x \cap \partial F \vert = i \}$$ 

Therefore  $l  =  l_0 + l_1 + l_2$ where $l_i = $ length of $g_i, i = 0,1,2$.
 
Each of the sets $g_i$ is a collection of intervals. If $A_1 = $ Area of $F$
outside $N_{\e}$ (cusp region) and $A_2 =$ area  of $F$ consisting of points
whose perpendicular geodesic to $g$ has length  less than $\e$ (narrow region) then
\begin{eqnarray*}
\mbox{Area}(F) & = &\mbox{Area}(F \cap N_{\e})  + A_1 \\
\mbox{Area}(F \cap N_{\e}) & = & 2l_0 \sinh(\e) + l_1 \sinh(\e) + A_2 \\
\Rightarrow  l \sinh(\e) &  = & 2l_0 \sinh(\e) + l_1 \sinh(\e) + A_1 + A_2 \\
\Rightarrow  l_2 \sinh(\e) &  = & l_0 \sinh(\e) + A
\end{eqnarray*}
where $A = A_1 + A_2$ (cusp plus narrow region).

\begin{figure}
\begin{center}~
\psfig{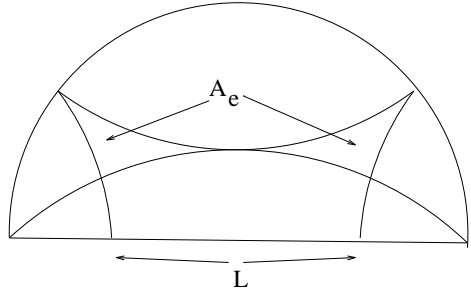}
\end{center}
\caption{$A_e$ is minimum \label{tangent}}
\end{figure}

If $e$ is one of the $n-2$ edges of $F$ not intersecting $g$, drop
perpendiculars to $g$ from the endpoints of $e$ on $S^1_{\infty}$ to the two
points $x,y$ on $g$. This forms a quadrilateral $Q_e$ with area $\pi$ and
base $[x,y]$. Let $A_e$ be the area of the region in $Q_e$ which is either a 
distance
greater than $\e$ from $g$ or whose perpendicular to $g$ has length less than
$\e$. $A_e$ is the contribution of $Q_e$ to $A$.
If $e$ doesn't intersect $N_{\e}$ then the $A_e$ is minimum when $e$ is 
tangent to $\partial N_{\e}$ (figure \ref{tangent}). If $L$ is the length of
the base of the quadrilateral tangent to $\partial N_{\e}$ then
\begin{eqnarray*}
& \sinh(L/2) & =  \frac{1}{\sinh(\e)} \\
\Rightarrow & A_e &  \geq  \pi - L \sinh(\e) = \pi - \frac{L}{\sinh(L/2)} \\
\Rightarrow & A_e & \geq  \pi - 2
\end{eqnarray*}

\begin{figure}
\begin{center}~
\psfig{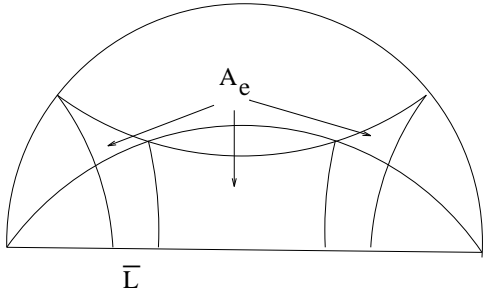}
\end{center}
\caption{General picture of $A_e$ \label{smallarea}}
\end{figure}

If $e$ intersects $\partial N_{\e}$ at points $x',y'$ then project down onto 
$g$ to get points $x_1,y_1$ (figure \ref{smallarea}). Now if $\o{L}$ is the 
length of either $[x,x_1]$ or $[y_1,y]$ (both have the same length) then 
\begin{eqnarray*}
\o{L} & \leq & L/2 \\
\Rightarrow A_e & = & \pi - 2\o{L} \sinh(\e) \\
\Rightarrow A_e & \geq & \pi - 2 \\
\Rightarrow A & \geq & (n-2)(\pi - 2)
\end{eqnarray*}

Since we have that $l_2 \sinh(\e) = l_0 \sinh(\e) + A$ then
\begin{eqnarray*}
l_2 & \geq & \frac{A}{\sinh(\e)} \geq \frac{(n-2)(\pi - 2)}{\sinh(\e)} \\
\Rightarrow l_2 & \geq & \frac{\pi -2}{\pi} \cdot l
\end{eqnarray*}
The subset $g_2$ consists of at most $(n-3)$ intervals. Therefore 
by the above, one of 
these intervals $[a,b]$ must have length $D$ satisfying
$$D \geq \frac{\pi -2}{\pi(n-3)} \cdot l$$
Consider the quadrilateral $Q$ formed by taking the two perpendicular arcs
$\c_a,\c_b$ to $g$, centered at $a,b$ respectively, of length $2\e$.
Choose the other sides $\a,\b$
to join the endpoints of $\c_a,\c_b$. There are two edges 
$e_3,e_4$ of $F$ separated by $g$, which enter and leave $Q$ through the
sides $\c_a,\c_b$. To show that $e_3,e_4$ are close note that they are closer
than $\a$ and $\b$ are. The two geodesics $\a,\b$ describe an ideal 
quadrilateral by extending $\a,\b$ to the sphere at infinity.
This quadrilateral has a diagonal of
length $D'$
which is longer than $D$ whose other diagonal has length $d'$. If 
$d$ is the distance between $e_3$ and $e_4$ then
$$\sinh(D/2) \leq \sinh(D'/2) = \frac{1}{\sinh(d'/2)}
\leq \frac{1}{\sinh(d/2)}$$

Therefore we have
\begin{equation}
\sinh(d/2) \cdot \sinh(\frac{(\pi - 2)l}{2\pi (n-3)}) \leq 1 
\label{longthin}
\end{equation}
 \eproof

This lemma shows that long geodesics must give thin polygons.

\subsection*{Thin Bent Plane Forces Intersection}
In the hyperbolic plane consider a geodesic $\c$ represented by a circle $C$ 
of radius $R$ intersecting the unit circle at right angles. 
The two geodesics  passing 
through the origin $o$ meeting $\c$ on $S^1_{\infty}$ form an angle of $\f$ 
at $o$ (visual angle). If $x$ is the Euclidean distance from $o$ to $\c$ and
$d$ is the hyperbolic distance then we have by simple planar geometry
\begin{eqnarray*}
& R & = \frac{1 - x^2}{2x} \\
& \tan(\f/2) & = R \\
& 2\tanh^{-1}(x) & = d \\
\Rightarrow & x & = \tanh(d/2) \\
\Rightarrow & R & = \frac{1}{\sinh(d)} \\
\Rightarrow & \tan(\f/2) & = \frac{1}{\sinh(d)}
\end{eqnarray*}

In $P_{\t}$ we have planes $F_1,F_2$ intersecting with dihedral angle $\f$ and
let $F_3,F_4$ be the two planes intersecting in edges $e_3,e_4$ respectively of
lemma \ref{thin}. We also have planes $F'_3,F'_4$ defined by geodesics
 $\a,\b$ which intersect $F_1,F_2$ at right angles.

\begin{figure}
\begin{center}~
\psfig{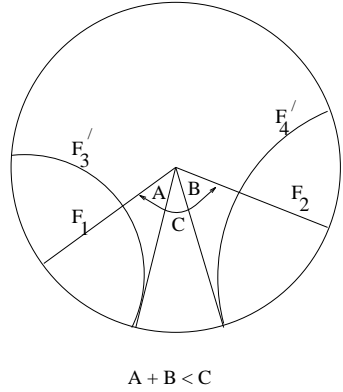}
\end{center}
\caption{Sideview of the planes \label{sideview}}
\end{figure}

If $F'_3,F'_4$ intersect then  $F_3,F_4$ must also intersect. Therefore
since planes $F_3,F_4$ cannot 
intersect then neither can $F'_3,F'_4$. Looking along the geodesic 
$g$,   $F'_3,F'_4$ are given by two geodesics $\c_3,\c_4$ and $F'_3,F'_4$   
by two geodesics
passing through the origin and perpendicular to $\c_3,\c_4$ 
respectively (figure \ref{sideview}).

Both $\c_3,\c_4$ have the same visual angle $\f'$. If they are not to 
intersect then
\begin{eqnarray*}
\f'/2 & \leq & \f/2 \\
\Rightarrow \tan(\f'/2) & \leq & \tan(\f/2) \\
\Rightarrow \frac{1}{\sinh(d'/2)} & \leq & \tan(\f/2)
\end{eqnarray*}
where $d'$ is the length of the short curve on $F$ joining $F'_3,F'_4$. By
lemma \ref{thin} we have
$$\sinh(\frac{(\pi - 2)l(\f)}{2\pi(n-3)}) \leq \tan(\f/2)$$
Therefore
\begin{equation}
l(\f) \leq \frac{2\pi(n-3)}{(\pi - 2)} \cdot \sinh^{-1} (tan(\f/2))
\label{bfn}
\end{equation}

Therefore we obtain a combinatorial bound on the volume 
increase by integrating the above equation. This integration is possible
as although $\sinh^{-1}(\tan(\f/2))$ tends to infinity at $\pi$ it is still
integrable. Therefore we have
\begin{eqnarray}
\Delta V & \leq & \frac{\pi(n-3)}{(\pi - 2)} \cdot \int_{0}^{\pi} \sinh^{-1}(\tan(\f/2)d\f = K \cdot (n-3)
\end{eqnarray}
where $K$ is some fixed constant.
We can derive this bound in general using the one-parameter family of 
approximations $X^t$ in the same manner as before.

\section{Non-existence of lower bound for volume increase}

If the volume increase under Dehn drilling is bounded below by a constant times
some positive power of the length of the geodesic drilled then we have that
$$\Delta V \geq c \cdot l^a$$ where $c,a$ are fixed positive constants and 
$l$ is the length of the geodesic drilled.
We have shown in the previous section, that the volume increase by drilling out
a geodesic belonging to an ideal $n$-gon of a basic polyhedron is bounded 
linearly in $n$.  Therefore to prove that no lower bound exists we need
to find a sequence of basic polyhedra each containing a face of bounded 
valence
and the set of lengths of geodesics perpendicular to two sides unbounded. Let 
$\{P_k\}_{k=0}^{\infty}$
be a sequence  of basic polyhedra and $N$ a number  such that
each $P_k$ has a marked face $F_k$ which is an $n_k$-gon, $n_k < N$. Also
let $l_k$ be the length of the longest geodesic $\a_k$ 
in $F_k$ perpendicular to two edges of $F_k$. Then 
the volume increase $\Delta V_k$, by 
drilling out $\a_k$ is
bounded above by some constant $K$ depending only on $N$. If 
the sequence $l_k$ is 
not bounded then by the lower bound for volume increase 
$\Delta V_k \geq c \cdot l_k^a$. Therefore $\Delta V_k$ is 
unbounded giving a contradiction.
We now construct the required family using circle packings.

\subsection*{Circle Patterns and Packings}
If $P$ is a hyperbolic polyhedron then in the upper half-space model 
each face $F$ of $P$ lies on a
hemisphere $H_F$ which intersects $S^2_{\infty}$ in a circle $C_F$. Furthermore
the angle at which two circles $C_{F_1},C_{F_2}$ intersect equals the dihedral 
angle  between the faces $F_{1},F_2$. Thus looking for a polyhedron of 
prescribed 
dihedral angles and combinatorics is equivalent to finding a
circle pattern with prescribed angles of intersection and 
combinatorial pattern. We will use the terminology circle packing to refer to
a collection of circles with intersections being tangential.

\begin{figure}
\begin{center}~
\psfig{file=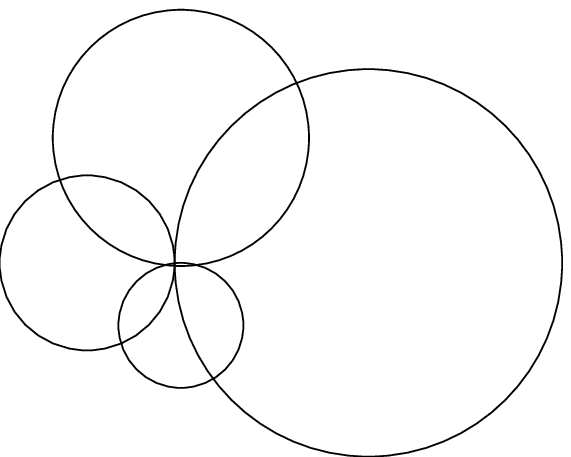,height=1.0in}
\end{center}
\caption{\label{4circ}}
\end{figure}

In the case of constructing a basic polyhedron we require that the circles meet
in groups of four with opposite circles tangent and neighboring circles 
meeting at right-angles (figure \ref{4circ}).

We start our circle pattern with a circle $C_0$ of radius 
one centered at the origin 
$O$. We further take circles $C_{1},C_{1}'$ of radius $r$ such that both are
perpendicular to $C_0$ and have centers on the positive, negative real
axis respectively. Thus there are unique circles $C_{2},C_{2}'$ which are 
perpendicular
to $C_0$, with centers on the positive, negative  imaginary  axis
respectively and are tangent to both $C_{1},C_{1}'$. This 
configuration of five circles depends only on our choice of $r$ and therefore 
we have 
a one-parameter family ${\mathcal{C}}_r$ of circle patterns. By simple 
planar hyperbolic geometry the radius of both  $C_{2},C_{2}'$ is $1/r$. 
This collection of five circles
describes an ideal quadrilateral on the hyperbolic plane $H_0$ associated with
circle $C_0$. Furthermore 
$r$ is a coordinate system for the shape of the quadrilateral with it becoming
narrow in the horizontal direction as $r$ tends to 
infinity and becoming narrow 
in the vertical direction $r$ tends to $0$.

\begin{figure}
\begin{center}~
\psfig{file=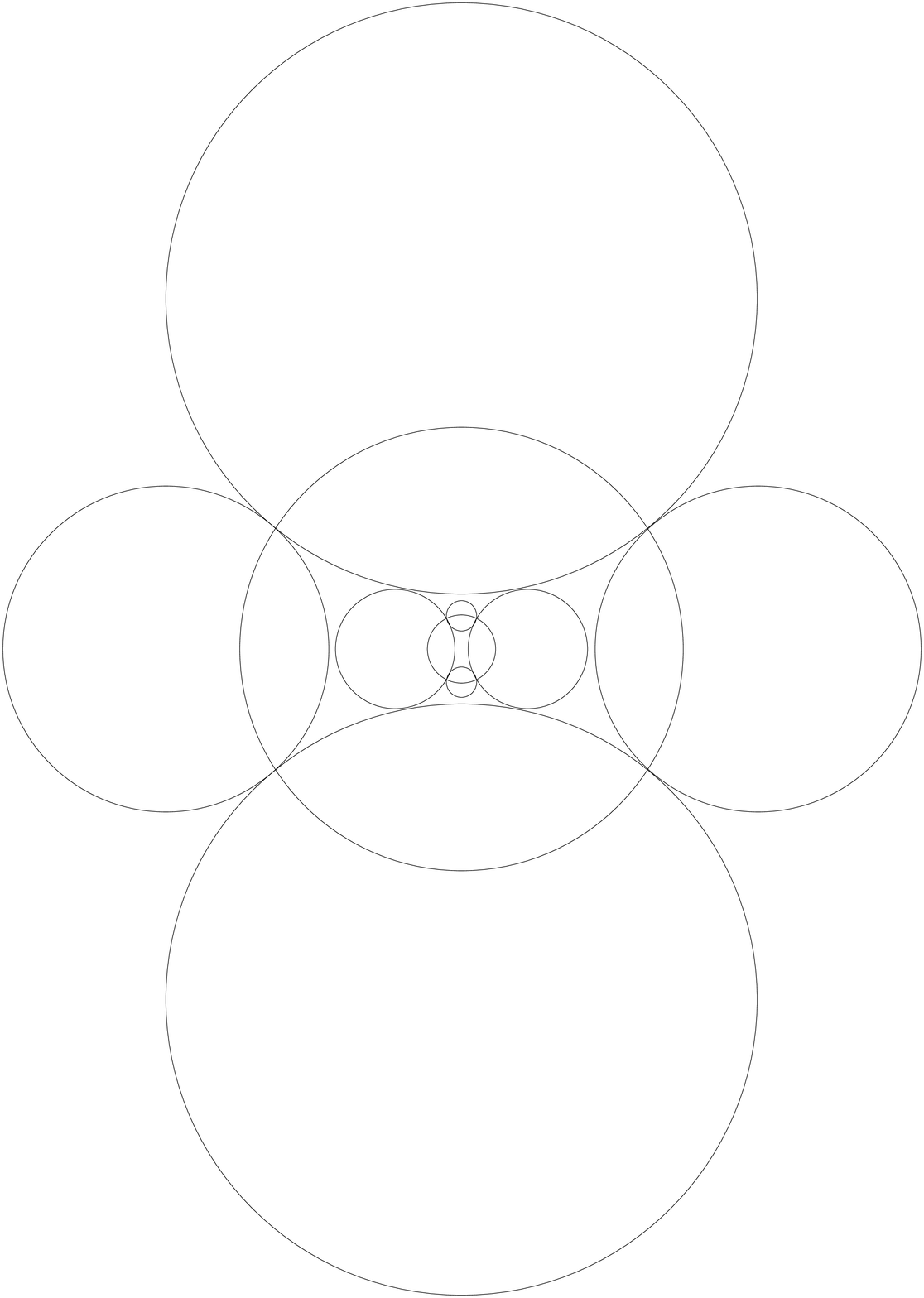,height=8.0in}
\end{center}
\caption{Circle pattern \label{cpacking}}
\end{figure}

To extend ${\mathcal{C}}_r$ to the circle pattern of a basic polyhedron 
we perform
a moebius transformation $M$ on ${\mathcal{C}}_r$. $M$ consists of a 
rotation of 
$\pi/2$ about the origin followed by an expansion by a real number $k_r$.
The factor $k_r$ is chosen so that  $M(C_{1})$ is tangent
to the original $C_1$. For $r$ large these four tangent points are the only 
intersections between ${\mathcal{C}}_r$ and $M({\mathcal{C}}_r)$
and so we  extend  circle pattern ${\mathcal{C}}_r$ by adding
$M({\mathcal{C}}_r)$ to form circle pattern 
$\o{{\mathcal{C}}_r} = {\mathcal{C}}_r  \cup  M({\mathcal{C}}_r)$ 
(figure \ref{cpacking}).

Now take the collection of hemispheres ${\mathcal{H}}_r$ associated with 
$\o{{\mathcal{C}}_r}$ and let $N_r$ be the region between the two hemispheres 
$H_{0},M(H_{0})$ cut out by ${\mathcal{H}}_r$. $N_r$ is an infinite volume 
hyperbolic orbifold with four
infinite volume funnels. Each of these funnels meet $S^2_{\infty}$ 
in a quadrilateral bounded by  circles which are tangent at points of 
intersection. It can be easily shown that all four of these quadrilaterals
are equivalent up to moebius transformation to a quadrilateral $Q_r$.

To form a basic polyhedron from $N_r$ we must truncate the funnels of $N_r$. 
We do this by first finding  a packing of $Q_r$ by a 
finite number of circles such that the interstices are all triangular. Then
having packed each of the four copies of $Q_r$,  we are left with a finite 
number of triangular interstices which we truncate by 
adding the unique circle 
passing through the three vertices of each  interstice. These  circles are
also perpendicular to the three sides of the 
interstices and therefore we have truncated the four infinite volume funnels of
$N_r$ to get a basic polyhedron $P_r$.

Unfortunately it is not always possible to find a finite packing of an 
arbitrary circular quadrilateral but using work of Robert Brooks(\cite{Br})
we will show that there is a sequence $\{r_{n}\}$ such that $Q_{r_{n}}$ is 
packable and $r_{n}$ monotonically tends to infinity.

\begin{figure}
\begin{center}~
\psfig{file=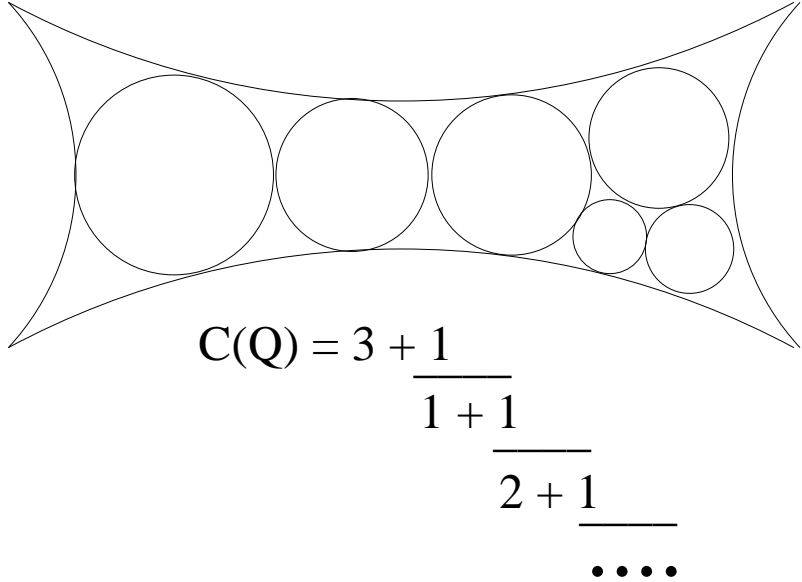,height=3.0in}
\end{center}
\caption{Greedy algorithm\label{cfrac}}
\end{figure}

We define the continued fraction $c(Q)$ of a quadrilateral $Q$ by using the 
greedy algorithm.
Assume that the four sides have been labelled left, right, top, bottom. 
We first 
pack $Q$ from left to right by starting at the left
side of $Q$ and  packing  in the largest circle tangent to the left, top
and bottom of $Q$. Now $Q$ contains a smaller quadrilateral $Q_1$ 
which we again
pack with the largest circle tangent to the left, top and bottom sides.
After $n_1$ circles have been added the remaining quadrilateral $Q_{n_1}$
is too narrow
to fit any more circles which are tangent to the left, top and bottom so  
we now pack  $Q_{n_1}$ from top to bottom with   
circles tangent to the top, left and right  sides. Again after $n_2$ circles
have been packed into $Q_{n_1}$ we are left with a quadrilateral 
$Q_{n_{1},n_{2}}$
too wide  to pack with a circle tangent to the top, left and right. The 
packing algorithm is defined recursively by now applying the algorithm to
the quadrilateral $Q_{n_{1},n_{2}}$ (figure \ref{cfrac}). To each 
quadrilateral we get a sequence of
integers $\{n_{i}\}$ and we define the function $c$ on the space of 
quadrilaterals to be the real number with these as its continued fraction 
expansion, that is

\begin{eqnarray*}
 c(Q) = n_{1} + \lcfrac{1}{n_{2} +
		 \lcfrac{1}{n_{3} +
		  \lcfrac{1}{n_{4} + \dotsb
		}}}      
\end{eqnarray*}

Note that the greedy algorithm gives a finite packing of $Q$ if and only if 
$c(Q)$ is rational. Now let  $\mathcal{S}$ be the space of all 
configurations of four circles up to moebius transformation 
which have  disjoint interiors and such that 
the complement of the interiors is two circular quadrilaterals. Then 
Robert Brooks showed (\cite{Br}) that the map 
$\o{c} : \mathcal{S} \rightarrow \Re^2$ defined
by mapping a configuration of circles to the pair of continued fractions of
its two circular quadrilateral regions, is not only continuous but a 
homeomorphism. This implies that  $c$ is continuous and therefore
$c(Q_r)$ is a continuous function of $r$. It can be shown that as
$r$ tends to infinity
the modulus of $Q_r$ tends to infinity. Therefore  $c(Q_r)$ tends to infinity
and by continuity there exists a sequence $\{r_{n}\}$ such that 
$c(Q_{r_n}) = n$. Therefore $Q_{r_{n}}$ can be packed with $n$ circles from 
left to right. 

\begin{figure}
\begin{center}~
\psfig{file=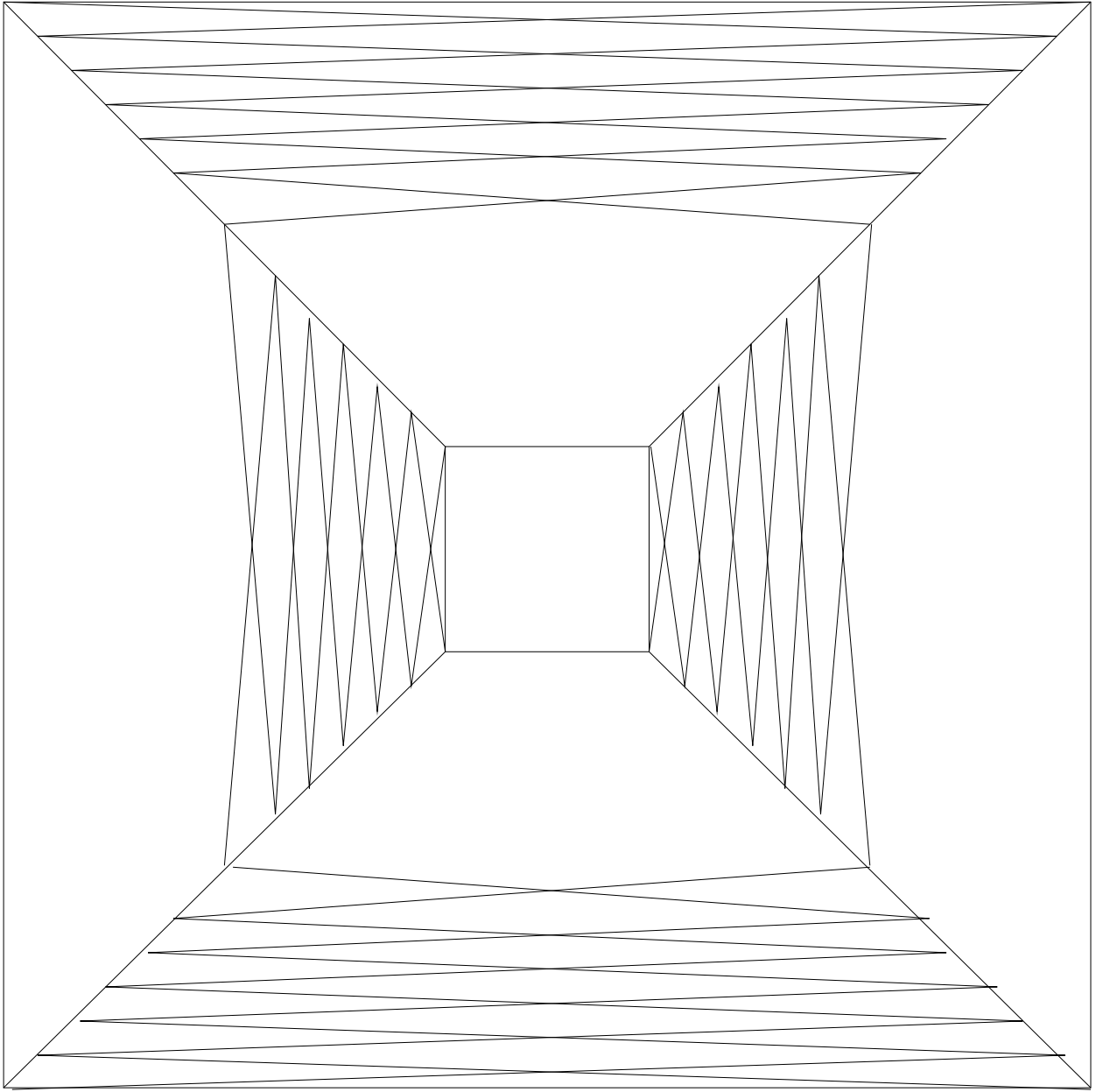,height=6.5in}
\end{center}
\caption{Polyhedron with thin ideal quadrilateral face\label{poly}}
\end{figure}

Now basic polyhedron $P_n$ is defined by truncating the funnels of 
$N_{r_{n}}$  as described above. Figure \ref{poly} is a combinatorial
picture of $P_n$. The polyhedra $P_n$ have the property of 
each containing an ideal
quadrilateral face $F_n$ whose modulus tends to infinity as a function
in $n$. Thus $F_n$ contains a geodesic $\a_n$, perpendicular to two sides
of $F_n$ whose length monotonically increases to infinity as a function of $n$.
Therefore this sequence of
polyhedra is a counter-example to the existence of a lower bound 
to volume increase under Dehn drilling. We get a similar counter-example for 
the manifold case by looking at the four-fold covers of 
these polyhedra when viewed as orbifolds.

\subsection*{Conclusion}
The two bounds we derive are very different. One is a
bound on the volume increase by the length of the geodesic drilled and the 
other is a purely combinatorial bound.
I believe that these techniques can be applied to a wider
class than I have done and further hope to prove
the conjecture at the outset by decomposing the manifold into polyhedra
and extending the results in this paper.

\end{document}